\documentclass[11pt, reqno]{amsart}

\usepackage{amsmath, amsthm, amscd, amsfonts, amssymb, graphicx, color}
\usepackage[bookmarksnumbered,plainpages]{hyperref}

\newtheorem{thm}{Theorem}[section]
\newtheorem{cor}[thm]{Corollary}
\newtheorem{lem}[thm]{Lemma}

\theoremstyle{definition}
\newtheorem{defn}[thm]{Definition}
\newtheorem{rem}[thm]{Remark}
\numberwithin{equation}{section}

\newcommand{\A}{\mathcal{A}}
\newcommand{\E}{E_{\epsilon}}
\newcommand{\Ed}{E_{\delta}}

\newcommand{\f}{\mathcal{F}}

\newcommand{\B}{\mathcal{B}}
\newcommand{\s}{S^{Lmc}}

\newcommand{\LM}{Lmc}

\newcommand{\p}{\widetilde{p}}

\title[]{$Lmc-$compactification of a semitopological semigroup as a space of e-ultrafilters }
\author[]{M. Akbari Tootkaboni,\\
 Department of Mathematics,
 Shahed university \\
 Tehran, Iran\\
 Email: akbari@shahed.ac.ir}
\begin{document}
\begin{abstract}
Let $S$ be a semitopological semigroup and $\mathcal{CB}(S)$
denotes the $C^*$-algebra of all bounded complex valued continuous
functions on $S$ with uniform norm. A function $f\in \mathcal{CB}(S)$ is left multiplicative \linebreak continuous
if and only if $\mathbf{T}_{\mu}f\in \mathcal{CB}(S)$ for all $\mu$ in the spectrum of $\mathcal{CB}(S)$,
where $\mathbf{T}_{\mu}f(s)=\mu(L_sf)$ and $L_sf(x)=f(sx)$ for each $s,x\in S$.
The collection of all left multiplicative continuous functions on $S$ is denoted by $Lmc(S)$.
In this paper, the $Lmc-$compactification of a semitopological semigroup S is reconstructed as a space of $e-$ultrafilters.
This construction is applied to obtain some algebraic properties of \linebreak $(\varepsilon ,\s)$, that $ \s $ is the
spectrum of $Lmc(S)$, for  semitopological semigroups $S$. It is shown that if  S is a locally compact
semitopological semigroup, then   $S^*=\s \setminus \varepsilon(S)$ is a left ideal
of $\s$ if and only if for each $x,y\in S$, there exists a compact zero set $A$ such that $x\in A$
and $\{t\in S:yt\in A\}$ is a compact set.
\end{abstract}
\maketitle

\section{\bf Introduction}

Let $X$ be a completely regular space, $\mathcal{C}(X,\mathbb{R})$ denotes all the
real \linebreak valued continuous function on $X$ and $\mathcal{CB}(X,\mathbb{R})$
denotes all the bounded real valued continuous function on $X$. The
correspondences between $z-$filters on $X$ and ideals in
$\mathcal{C}(X,\mathbb{R})$ that have been established in \cite{Gil} are
powerful tools in the study of $\mathcal{C}(X,\mathbb{R})$. These
correspondences, which also occur in a rudimentary form in
$\mathcal{CB}(X,\mathbb{R})$, are inconsequential, as many theorems of  \cite{Gil} become false if $\mathcal{C}(X,\mathbb{R})$ is replaced
by $\mathcal{CB}(X,\mathbb{R})$. However, there is \linebreak another correspondence between a certain class of $z-$filters on $X$ and ideals in
$\mathcal{CB}(X,\mathbb{R})$ that leads to quite analogous theorems to those
for $\mathcal{C}(X,\mathbb{R})$. The requisite information is outlined in
\cite[2L]{Gil}.

In Section 2,  some familiarity with semigroup compactification
 and \linebreak $Lmc-$compactification will be presented. Also, this Section consists of an
introduction to $z-$filters and an elementary external
construction of $\LM-$compactification  as a space of $z-$filters. Moreover, in this section  $e-$filters and $e-$ideals will be defined (this has been adopted
from \cite[2L]{Gil}.

In Section 3, $\LM-$compactification will be reconstructed as a space of  $e-$ ultrafilters with a suitable topology,
 also  on $e-$ultrafilters a binary operation will be defined.

Section 4 will be about some theorems from \cite{hindbook}
 about the properties of $\beta S$ which are extended to some properties on $\s$, for semitopological semigroup S.

\section{\bf Preliminary}

Let $S$ be a semitopological semigroup (i.e. for each $s\in S$,
$\lambda_{s} :S\rightarrow S$  and $r_{s} :S\rightarrow S$ are
continuous, where for each $x\in S$, $\lambda_{s}(x)=sx$ and
$r_{s}(x)=xs$) with a Hausdorff topology, $\mathcal{CB}(S)$
denotes the $C^*$-algebra of all bounded complex valued continuous
functions on $S$ with uniform norm, and $\mathcal{C}(S)$
denotes the algebra of all complex valued continuous
functions on $S$. A  \mbox{semigroup} compactification of $S$ is a pair
$(\psi,X)$, where $X$ is a compact, Hausdorff, right topological
semigroup (i.e. for all $x\in X$, $r_x$ is continuous) and
$\psi:S\rightarrow X$ is continuous homomorphism with dense image
such that, for all $s\in S$, the mapping
$x\mapsto\psi(s)x:X\rightarrow X$ is \mbox{continuous}, (see Definition
3.1.1 in \cite{Analyson}). Let $\mathcal{F}$ be a $C^*$-subalgebra
 of $\mathcal{CB}(S)$ containing the constant functions, then the set of all \mbox{multiplicative} means of
 $\f$ (the spectrum of $\f$), denoted by $S^{\f}$, equipped with the
 Gelfand topology, is a compact Hausdorff topological space. Let $R_sf=f\circ r_{s}\in\f$ and $L_sf=f\circ \lambda_{s}\in\f$ for all
$s\in S$ and $f\in\f$, and the \mbox{function} $s\mapsto(T_\mu
f(s))=\mu(L_s f)$ is in $\f$ for all $f\in\f$ and $\mu\in
S^{\f}$ then $S^{\f}$ under the multiplication
$\mu\nu=\mu\circ T_\nu$ ($\mu,\nu\in S^{\mathcal{F}})$, furnished
with the Gelfand topology, makes $(\varepsilon, S^{\f})$ a
semigroup compactification (called the $\f$-compactification) of
$S$, where $\varepsilon :S\rightarrow S^{\mathcal{F}}$ is the
evaluation mapping. Also,
$\varepsilon^{*}:\mathcal{C}(S^{\mathcal{F}})\rightarrow\mathcal{F}$
is  isometric isomorphism and \mbox{
$\widehat{f}=(\varepsilon^{\ast})^{-1}(f)\in
\mathcal{C}(S^{\mathcal{F}})$} for $f\in\f$ is given by
$\widehat{f}(\mu)=\mu(f)$ for all $\mu\in S^{\f}$, (For more
detail see section 2 in \cite{Analyson}).

Let $\f=\mathcal{CB}(S)$, then $\beta S=S^{\mathcal{CB}(S)}$ is
the Stone-$\check{C}$ech compactification of $S$, where $S$ is a completely regular space.

A function $f\in \mathcal{CB}(S)$ is left multiplicative continuous
if and only if $\mathbf{T}_{\mu}f\in \mathcal{CB}(S)$ for all
$\mu\in\beta S=S^{\mathcal{CB}(S)}$. Therefore
$$Lmc(S)=\bigcap\{\mathbf{T}^{-1}_{\mu}(\mathcal{CB}(S)):\mu\in\beta S\}$$
is defined and $(\varepsilon,S^{Lmc})$ is the
universal semigroup compactification of \linebreak $S$ (Definition 4.5.1 and Theorem
4.5.2 in \cite{Analyson}). In general, $S$ can not be  embedded in $\s$. In fact, as it was shown in \cite{hin.mil} there is
a completely regular Hausdorff semitopological semigroup $S$ such that the
continuous \mbox{homomorphism} $\varepsilon$ from $S$ to its $Lmc-$compactification, is neither \mbox{one-to-one} nor
open as a mapping to $\varepsilon(S)$.

The $\mathcal{LUC}-$compactification is the spectrum of the $C^*-$algebra consisting of all left uniformly continuous
functions on semitopological semigroup $S$; a function $f:S\rightarrow \mathbb{C}$ is left uniformly continuous if
$s\mapsto L_sf$ is continuous map from $S$ to the space of bounded continuous functions on $S$ with the uniform norm. Let $G$ be a locally compact Hausdorff
topological group, by Theorem 5.7 of chapter 4  in \cite{Analyson} implies that $\LM(G)=\mathcal{LUC}(G)$.
Also the evaluation map $G\rightarrow G^\mathcal{LUC}$ is open, ( see \cite{Analyson}).

Now, some prerequisite material from \cite{Akbari} are quoted for the
description of $(\varepsilon,S^{Lmc})$ in terms of $z-$filters. For
$f\in\LM(S)$, $Z(f)=f^{-1}(\{0\})$ is called zero set,
and the collection of all zero sets is denoted by $Z(\LM(S))$. For an
extensive account of ultrafilters, the readers may refer to
\cite{Ultra},
 \cite{Gil} and \cite{hindbook}.

\begin{defn}
$\mathcal{A}\subseteq Z(\LM(S))$ is called a $z$-filter on $\LM(S)$ (or for simplicity $z-$filter) if,
\\
(i) $\emptyset\notin \mathcal{A}$ and $S\in\A,$\\
(ii) if $A,B\in \mathcal{A}$, then $A\bigcap B\in\mathcal{A}$,
\\
(iii) if $A\in \mathcal{A}$, $B\in Z(\LM(S))$ and $A\subseteq
B$ then $B\in \mathcal{A}$.
\end{defn}
Because of $(iii)$, $(ii)$ may be replaced  by,
$(ii)'$ If $A,B\in \A$, then $A\cap B$ contains a member of
$\A$.

A $z-$ultrafilter is a  $z-$filter which  is not properly  contained
 in any other $z-$filter. The collection of all $z-$ultrafilters is denoted by $\mathcal{Z}(S)$. For $x\in S$,  $\widehat{x}=\{Z(f):f\in\LM(S),\,\,f(x)=0\}$
 is a $z-$ultrafilter. The $z-$filter $\f$ is said to converge to the limit $\mu\in S^{Lmc}$ if every neighborhood
 of $\mu$ contains a member of $\f$. The collection of all $z-$ultrafilters on $\LM(S)$ converge to $\mu\in S^{Lmc}$ is
 denoted by $[\mu]$. Let $\mathcal{Q}= \{\widetilde{p}:\p=\cap[\mu]\}$ and define $\widetilde{A}=\{\p:A\in \p\}$ for $A\subseteq S$. Let
$\mathcal{Q}$ be equipped with the topology whose basis is
 $ \{(\widetilde{A})^{c}: A\in Z (\LM(S))\} $, and define $\bigcap [\mu]*\bigcap [\nu]=\bigcap [\mu\nu]$.  Then $ (\mathcal{Q}, *) $ is  a (Hausdorff)
 compact right topological semigroup and $ \varphi : \s \rightarrow \mathcal{Q}  $ defined by $ \varphi( \mu)= \cap [\mu]=\widetilde{p} $,
where $\bigcap_{A\in p}\overline{A}=\{ \mu \}$,  is   topologically isomorphism. So $ \widetilde{A} $ is equal  to
$ cl_{ S^{Lmc}}A $  and we denote it by $ \overline{A} $, also for simplicity we use $x$ replace $\widehat{x}$.
 The operation $ ``\cdot" $ on $  S $, extends uniquely to $ ``*" $ on $ \mathcal{Q} $.
 For  more discussion and  details see
 \cite{Akbari}.

\begin{rem}
 If $p,q\in \mathcal{Z}(S)$, then the following statements hold.
\\
(i) If $E\subseteq Z(\LM(S))$ has the finite intersection property, then $E$ is
 \mbox{contained} in a $z-$ultrafilter.
\\
(ii) If $B\in Z(\LM(S))$ and for all $A\in p$, $A\cap B\neq\emptyset$ then
$B\in p$.
\\
(iii) If $A, B\in Z(\LM(S))$ such that $A\cup B\in p$ then $A\in p$ or
$B\in p$.
\\
(iv) Let $p$ and $q$ be distinct $z-$ultrafilters, then there exist $A\in p$ and $B\in q$ such that
$A\cap B=\emptyset$.
\\
(v) Let $p$ be a $z-$ultrafilter, then there exists $\mu\in\s$ such that
$$\bigcap_{A\in p}\overline{\varepsilon(A)}=\{\mu\}.$$
(For (i), (ii), (iii) and (iv) see Lemma 2.2  and Lemma 2.3 in \cite{Akbari}. For (v) see Lemma 2.6 in \cite{Akbari}).
\end{rem}

In this paper, $\mathbb{R}$ denotes the topological group formed by the real numbers
under addition. Also define $Ker(\mu)=\{f\in\LM(S):\mu(f)=0\}$ for $\mu\in \LM(S)^*$. By
Theorem 11.5 in \cite{Rud}, $M$ is a maximal ideal of $\LM(S)$ if and only if there
is $\mu\in\s$ such that $Ker(\mu)=M$.

\begin{lem}
Let $S$ be a Hausdorff semitopological semigroup.
\\
(1) If $f\in\LM(S)$ is  real valued, then \mbox{
$f^+=max\{f,0\}\in\LM(S)$} and $f^-=-min\{f,0\}\in\LM(S)$.
\\
 (2) Let $f\in \LM(S)$, then $Re(f)\in \LM(S),\ Im(f)\in \LM(S)\text{ and }\linebreak |f|\in\LM(S) $.
\\
 (3) If $f \text{ and } g $ are real valued functions in $\LM(S)$, then $$(f\vee g)(x) =max\{f(x),g(x)\}\in \LM(S),$$ and  $$(f\wedge
g)(x)=min\{f(x),g(x)\}\in \LM(S).$$
\\
 (4) Let $f\in\LM(S)$
and there exists $c>0$ such that $c<|f(x)|$ for each $x\in S$. Then
$\frac{1}{f}\in\LM(S)$.
\\
(5) Let $M$ be a maximal ideal and $f\in M$, then $\overline{f}\in M.$
\end{lem}

{\flushleft\bf Proof.} For (1), (2) and (3), since $f\mapsto \widehat{f}:\LM(S)\rightarrow \mathcal{C}(S^{\LM})$ is
an isometrical isomorphism and $|\widehat{f}|\in \mathcal{C}(S^\LM)$
 for each $f\in\LM(S)$, so we have
 \begin{eqnarray*}
|\widehat{f}|(\varepsilon(x))&=&|\widehat{f}(\varepsilon(x))|\\
&=&|\varepsilon(x)(f)|\\
&=&|f(x)|\\
&=&|f|(x)
\end{eqnarray*}
for each $x\in S$. Thus, $|\widehat{f}|=\widehat{|f|}$ for each $f\in\LM(S)$ and so
$|f|\in\LM(S)$ for each $f\in\LM(S)$.\\
Now let $f$ and $g$ be real valued functions, so $$f\vee g(x)=\frac{|f-g|(x)}{2}+\frac{(f+g)(x)}{2}\in\LM(S).$$
In a similar way $f\wedge g,\,f^+ \text{ and } f^-$ are in $ \LM(S) $. Pick $f\in\LM(S)$, since $\LM(S)$ is conjugate closed subalgebra
so $Re(f)=\frac{f+\overline{f}}{2}\in\LM(S)$
and $Im(f)=\frac{f-\overline{f}}{2i}\in\LM(S)$.

For (4), let $f\in\LM(S)$
and there exists a $c>0$ such that $c<|f(x)|$ for each $x\in S$. So $|\widehat{f}|(\mu)\geq c$ for each $\mu\in S^{\LM}$, this implies
that $\widehat{\frac{1}{f}}=\frac{1}{\widehat{f}}\in \mathcal{C}(\s)$. Therefore, $\frac{1}{f}\in\LM(S)$.

For (5), let $M$ be a maximal ideal in $\LM(S)$, so there exists a \linebreak $\mu\in\s$ such that $M=Ker(\mu)=\{f\in\LM(S):\mu(f)=0\}$. Now let $f\in M$, so
$\mu(f)=\mu(Re(f))+i\mu(Im(f))=0$. This implies that $\mu(Re(f))=\mu(Im(f))=0$ and so $\mu(\overline{f})=0$. Thus, $\overline{f}\in M$.
$\hfill\square$

\vspace{.4cm}
For $f\in\LM (S)$ and $\epsilon >0$, we define
$E_{\epsilon}(f)=\{x\in S:|f(x)|\leq\epsilon\}$.
Every such
set is a zero set. Conversely, every zero set is of this form,
$Z(g)=\E (\epsilon +|g|).$
For $I\subseteq \LM (S)$, we write $E(I)=\{\E (f):f\in
I,\epsilon
>0\}$, i.e. $E(I)=\bigcup_{\epsilon>0}\E (I).$
 Finally, for any
family $\A$ of zero sets, we define  $$E^{-}(\A)=\{f\in \LM (S):
\E (f)\in\A \mbox{ for all }\epsilon >0\},$$ that is,
$E^{-}(\A)=\bigcap_{\epsilon>0}E^{\leftarrow}_{\epsilon}(\A),$ where $$E_\epsilon^{\leftarrow}(\A)=\{f\in\LM(S):\E(f)\in\A\}.$$

\begin{lem}
For any family $\A$ of zero sets,
$$E(E^{-}(\A))=\bigcup_{\epsilon>0}\{\E (f):f\in\LM(S),\,\,E_{\delta}(f)\in\A
\mbox{ for all }\delta >0\}\subseteq\A.$$   The inclusion may be
proper,  when $\A$ is a $z-$filter.
\end{lem}

{\flushleft\bf Proof.} Let $\A$ be a family of zero sets, so
$$E^{-}(\A)=\{f\in \LM (S):
\E (f)\in\A \mbox{ for all }\epsilon >0\},$$
and thus,
\begin{eqnarray*}
E(E^{-}(\A))&=&\{\E(f):f\in E^{-}(\A),\epsilon >0\}\\
&=&\bigcup_{\epsilon>0}\{\E(f):f\in E^{-}(\A)\}\\
&=&\bigcup_{\epsilon>0}\{\E(f):E_\delta(f)\in\A\mbox{ for all }\delta>0\}\\
&\subseteq&\A.
\end{eqnarray*}
Finally, let $M_0=\{f\in\LM((\mathbb{R},+)):f(0)=0\}$, then $M_0$ is a maximal ideal in $\LM((\mathbb{R},+))$ and $\A=\{Z(f):f\in M_0\}$
is a z-filter. Let $g(x)=|x|\wedge 1$,
then $g\in M_0$ and so $\{0\}=Z(g)\in\A$. Since
\begin{eqnarray*}
E(E^{-}(\A))&=&\bigcup_{\epsilon>0}\{\E(f):E_\delta(f)\in\A\mbox{ for all }\delta>0\}\\
&\subseteq&\A,
\end{eqnarray*}
pick $f\in\LM((\mathbb{R},+))$ such that $\E(f)\in E(E^-(\A))$ each $\epsilon>0$. Since $f$ is continuous
so for each $\epsilon>0$ there exists $\delta>0$ such that $f((-\delta,\delta))\subseteq (-\epsilon,\epsilon)$,
therefore $(-\delta,\delta)\subseteq\E(f)$. This implies that $E(E^-(\A))$ is a collection of uncountable zero sets.
But $\{0\}\in\A$ is finite and so $\{0\}\notin E(E^{-}(\A))$. Therefore
$E(E^{-}(\A))\neq\A$.
$\hfill\square$

\vspace{.2cm}
\begin{defn}
 Let $\A$ be a $z-$filter. Then $\A$ is called an $e-$filter if
$E(E^{-}(\A))=\A$.
\end{defn}

Hence, $\A$ is an $e-$filter if and only if, whenever $Z\in\A$,
there exist \linebreak $f\in\LM(S)$ and $\epsilon>0$ such that
$Z=E_{\epsilon}(f)$ and $E_{\delta}(f)\in\A$ for every
$\delta>0$.

\begin{lem}
Let $I$ be a subset of $\LM(S)$. Then, $$I\subseteq
E^{-}(E(I))=\{f\in\LM (S):E_{\epsilon}(f)\in E(I) \mbox{ for
all }\epsilon>0\}.$$ The inclusion may be proper, even when $I$ is
an ideal.
\end{lem}

{\flushleft\bf Proof.} By Definition $$I\subseteq
E^{-}(E(I))=\{f\in\LM (S):E_{\epsilon}(f)\in E(I) \mbox{ for
all }\epsilon>0\}.$$
Finally, let $I$ be the ideal of all functions
in $\LM((\mathbb{R},+))$ that vanish on a neighborhood of $0$. Pick $g(x)=|x|\wedge 1$ in $\LM((\mathbb{R},+))$ that vanishes precisely
at $0$. Since for each $\epsilon>0$, $\E(g)=E_{\frac{\epsilon}{2}}(g\vee\frac{\epsilon}{2}-\frac{\epsilon}{2})$
and $g\vee\frac{\epsilon}{2}-\frac{\epsilon}{2}\in I$, then $\E(g)\in E(I)$ for each $\epsilon>0$, and so
$g\in E^{-}(E(I))$ but $g\notin I$. This completes the proof.$\hfill\square$

\begin{defn} Let $I$ be an ideal of $\LM(S)$. $I$ is called an $e-$ideal if
$E^{-}(E(I))=I$.
\end{defn}

Hence, $I$ is an $e-$ideal if and only if, whenever
$E_{\epsilon}(f)\in E(I)$ for all $\epsilon>0$, then $f\in I$.

\begin{lem}
The following statements hold.
\\
 (1) The intersection of $e-$ideals
is an $e-$ideal.
\\
(2) If $I$ is an ideal in $\LM (S)$, then $E(I)$ is an $e-$ filter.
\\
(3) If $\A$ is any $z-$filter, then $E^{-}(\A)$ is an $e-$ideal in $\LM
(S)$.
\\
 (4) $I\subseteq J\subseteq \LM(S)$ implies $E(I)\subseteq E(J)$, and
$\A\subseteq\B\subseteq Z(\LM(S))$ implies $E^{-}(\A)\subseteq
E^{-}(\B)$.
\\
(5) If $J$
is an $e-$ideal, then $I\subseteq J$ if and only if
$E(I)\subseteq E(J)$. If $\A$ is an $e-$filter then $\A\subseteq \B$
 if and only if $E^{-}(\A)\subseteq E^{-}(\B)$.
\\
(6) If $\A$ is any $e-$filter, then $E^{-}(\A)$ is an $e-$ideal. Let
$I$ be an ideal in $\LM(S)$, then $E^{-}(E(I))$ is the smallest
$e-$ideal containing $I$. In particular, every maximal ideal in
$\LM(S)$ is an $e-$ideal.
\\
 (7) For any $z-$filter $\A$,
$E(E^{-}(\A))$ is the largest $e-$filter contained in $\A$.
\end{lem}

{\flushleft\bf Proof.} $(1)$ Suppose that  $\{I_\alpha\}$ is a collection of e-ideals and $I=\bigcap_\alpha I_\alpha$. Let
$\E(f)\in E(I)$ for each $\epsilon>0$, then $\E(f)\in E(I_\alpha)$ for each $\epsilon>0$ so,
$f\in I_\alpha$ for each $\alpha$. This implies $f\in I$.

$(2)$ Let $\E(f)=\emptyset$ for some $\epsilon>0$ and $f\in I$,
then $\epsilon \leq|f(x)|\leq M$ for some $M>0$ and for each $x\in S$. So $\frac{1}{f}\in \LM(S)$ and,
$1=f\frac{1}{f}\in I$. This is a contradiction and so $\emptyset\notin
E(I)$.

Let $f'\in\LM(S)$, $f\in I$ be a nonnegative function and
$\E (f)\subseteq Z(f')$, then
$g(x)=|f'(x)|+\frac{\epsilon}{\epsilon\vee |f(x)|}\in
\LM(S)$. Now
$$|f(x)|g(x)=|f'(x)f(x)|+\frac{\epsilon
|f(x)|}{\epsilon\vee |f(x)|},$$
 so for $x\in Z(f')$ implies that
$|f(x)g(x)|=\frac{\epsilon |f(x)|}{\epsilon\vee
|f(x)|}\leq\epsilon$. Hence we have $Z(f')\subseteq \E(fg)$. If
$x\in \E(fg)$, then $$|f'(x)f(x)|\leq |f'(x)f(x)|+\frac{\epsilon
|f(x)|}{\epsilon\vee |f(x)|}=|f(x)g(x)|\leq\epsilon$$ and if
$x\notin Z(f')$ then $\epsilon <|f(x)|$ and
$|g(x)f(x)|>\epsilon$. Therefore this implies $\E(fg)\subseteq Z(f')$,
and so $\E(fg)=Z(f')$.

Suppose that  $\E(f),\Ed(g)\in E(I)$ for some $f,g\in I$ and
$\epsilon,\delta>0$. Let $\gamma=\epsilon\wedge\delta\wedge
\frac{1}{2}$ then,
$$E_{\frac{\gamma^{2}}{4}}(f\overline{f}+g\overline{g})\subseteq E_{\gamma}
(f)\bigcap E_{\gamma}(g)\subseteq E_{\epsilon}(f)\bigcap E_{\delta}(g),$$
thus $\E(f)\cap\Ed(g)\in E(I)$.

Now let $Z\in E(I)$, so there exists a $f\in I$ such that $Z=\E(f)$ for some $\epsilon>0$. By Definition of $E(I)$,
$E_\delta(f)\in E(I)$ for each $\delta>0$, so $E(I)$ is an \linebreak $e-$ filter.

$(3)$ Let $f,g\in E^{-}(\A)$. Since $E_{\epsilon/2}(|f|)\cap
E_{\epsilon/2}(|g|)\subseteq \E(|f-g|) $, therefore $\E(f-g)\in\A$
for each $\epsilon>0$. Thus, $f-g\in E^{-}(\A)$.
Let $f\in E^{-}(\A)$, $g\in\LM(S)$ and $M=\|g\|+1$. Hence,
$E_{\frac{\epsilon}{M}}(f)\subseteq \E(fg)$ and $fg\in
E^{-}(\A)$. Now let $\E(f)\in E^{-}(\A)$ for each $\epsilon>0$.
Definition of $E^{-}(\A)$ implies that $f\in E^{-}(\A)$. Thus,
$E^{-}(\A)$ is an e-ideal.

$(4)$ It can easily be  checked.

$(5)$ It is obvious that if $I\subseteq J$ then $E(I)\subseteq E(J)$
by (4).

Conversely. If $f\in I$ then $\E(f)\in E(I)$ for each
$\epsilon>0$, and so $\E(f)\in E(J)$. Since $J$ is an
$e-$filter, so $f\in J$.
If $\A\subseteq \B$, then $E^{-}(\A)\subseteq E^{-}(\B)$. Since $\A$ is an $e-$filter,
then $\A=E(E^{-}(\A))\subseteq E(E^{-}(\B))\subseteq \B$.

$(6)$ Let $I=E^{-}(\A)=\{f\in \LM(S): \forall
\epsilon>0,\,\,\E(f) \in\A\}$, thus, $\A$ is an $e-$filter, and
$\A=E(E^{-}(\A))=E(I)$. This implies $I=E^{-}(\A)=E^{-}(E(I))$ and so
$I$ is an $e-$ideal.
  Let  $I\subseteq\LM(S)$ be an
ideal, then $J=E^{-}(E(I))$ is an  $e-$ ideal (by (3) and (4)), so $I\subseteq J$. Let $I\subseteq K\subseteq J$ and
$K$ be an $e-$ ideal, then $E(I)\subseteq E(K)\subseteq E(J)=E(I)$
and   $E(K)=E(I)$. Thus $J=E^{-}(E(I))=E^{-}(E(K))=K$, and this implies that $J$ is the smallest  $e-$ideal containing I.

Finally, every maximal ideal in $\LM(S)$ is an $e-$ideal. For this,
 let $M$ is a maximal ideal in $\LM(S)$. Then, $E^{-}
 (E(M))$ is an $e-$ideal, $M\subseteq E^{-}(E(M))$ and $M$ is maximal so,
 $M=E^{-}(E(M))$.

 $(7)$ Let $\A$ be a $z-$filter, then $E^{-}(\A)$ is an ideal in $\LM(S)$, so $E(E^{-}(\A))$
 is an $e-$filter and  $\B=E(E^{-}(\A))\subseteq
 \A$. Now let $\mathcal{U}$ be an $e-$filter such that $\B\subseteq
 \mathcal{U}\subseteq\A$, then $E^{-}(\mathcal{U})=E^{-}(\A)$.
 Hence, $\B$  is an $e-$ filter subset of $\A$.
 $\hfill\square$
\vspace{.3cm}

A maximal $e-$filter is called an $e-$ultrafilter. Zorn's Lemma implies that every $e-$filter
is contained in an $e-$ultrafilter.
 Because, if $\mathcal{Y}$ is a chain of $e-$filters, then also  is a chain of $z-$filters and
so $\cup \mathcal{Y}$ is a $z-$filter. It is sufficient to show $\cup\mathcal{Y}$ is an $e-$filter.
Let $Z\in\cup\mathcal{Y}$, then  there exists a $Y\in\mathcal{Y}$ such that $Z\in Y$. Since $Y$ is an $e-$ideal, so
 there exist $f\in\LM(S)$ and $\varepsilon>0$ such that $Z=\E(f)$ and $\{E_\delta(f):\delta>0\}\subseteq Y.$
Thus there exist $f\in\LM(S)$ and $\epsilon>0$ such that $Z=\E(f)$ and
$\{E_\delta(f):\delta>0\}\subseteq \cup\mathcal{Y}$. Therefore $\cup \mathcal{Y}$ is an $e-$filter.
\begin{thm}
If $M$ is a maximal ideal in $\LM(S)$  then $E(M)$ is an 
$e-$ ultrafilter, and if $\A$ is an $e-$ultrafilter then
$E^{-}(\A)$ is a maximal ideal in $\LM(S)$.
\end{thm}
{\flushleft\bf Proof.} Let $M$ be a maximal ideal, so, $E(M)$ is an
$e-$filter (Lemma 2.8(2)). Suppose that there exists an $e-$filter
$\mathcal{U}$ such that $E(M)\subseteq \mathcal{U}$, then \linebreak
$M=E^{-}(E(M))\subseteq E^{-}(\mathcal{U})$ and so,
$E(M)=E(E^{-}(\mathcal{U}))=\mathcal{U}$, by Lemma 2.8(7). Thus, $E(M)$ is an
$e-$ultrafilter.

 Now let $\mathcal{E}$ be an $e-$ultrafilter,
then $E^{-}(\mathcal{E})$ is an ideal in $\LM(S)$(Lemma 2.8(3)).
Let $J$ be a maximal ideal such that $E^{-}(\mathcal{E})\subseteq
J$, $J$ is an $e-$ideal and so, $
E(E^{-}(\mathcal{E}))\subseteq E(J)$. Since $\mathcal{E}$ is an
$e-$ultrafilter so, $\mathcal{E}=E(E^{-}(\mathcal{E}))$ and
 $E^{-}(\mathcal{E})=E^{-}(E(J))=J$. This implies that
$E^{-}(\mathcal{E})$ is maximal.$\hfill\square$
\vspace{.3cm}

The
correspondence $M\mapsto E(M)$ is one to one from the set of all
maximal ideals in $\LM(S)$ onto the set of all $e-$ultrafilters.
\begin{thm}
The following property characterizes an ideal $M$ in $\LM(S)$ as a
maximal ideal: given $f\in\LM(S)$,  if  $\E(f)$ meets every member of
$E(M)$ for each $\epsilon>0$, then $f\in M$.
\end{thm}
{\flushleft\bf Proof.} Let $M$ be a maximal ideal and $f\in\LM(S)$. Let $\E(f)$ meets every member of
$E(M)$ for each $\epsilon>0$. So $E(M)\cup\{\E(f):\epsilon>0\}$ has the finite intersection property,
and so there exists a $z-$ultrafilter $\A$ containing it. By Lemma 2.8 and Theorem 2.9,
$$M=E^{-}(\A)=\{g\in\LM(S):\E(g)\in\A\mbox{ for each }\epsilon>0\}.$$
This implies that $f\in M$.

Now let $M$ be an ideal in $\LM(S)$ with the following property: given $f\in\LM(S)$,  if  $\E(f)$ meets every member of
$E(M)$ for each $\epsilon>0$, then $f\in M$. We show that $M$ is a maximal ideal. Let
$f\in\LM(S)\setminus M$ and so some $\E(f)$ fails to meet some member of $E(M)$. Therefore there exist
$g\in M$ and $\delta>0$ such that $\E(f)\cap E_\delta(g)=\emptyset$. Pick $\gamma=min\{\delta^2,\epsilon^2,1\}$,
then $E_\gamma(f\overline{f}+g\overline{g})\subseteq E_{\epsilon}(f)\cap E_{\delta}(g)$, so
$f\overline{f}+g\overline{g}$ is invertible and generated ideal by $M\cup\{f\}$ is equal with $\LM(S)$. This implies $M$
is a maximal ideal.
$\hfill\square$

\vspace{.3cm}

Let $\A$ and $\B$ be $z-$ultrafilters. It is said that $\A \sim \B$ if and
only if $E(E^{-}(\A))=E(E^{-}(\B))$. It is obvious that $\sim$ is an
equivalence relation. The equivalence class of
$\A\in \mathcal{Z}(S)$ is denoted by $[\A]$.

\begin{lem}
Let $\A$ be a $z-$ultrafilter, then\\
(a) Let $Z(f)\in \A$ for some $\LM(S)$, then $f\in E^{-}(\A)$.\\
(b) $E^{-}(\A)$ is a maximal ideal.\\
(c) $E(E^{-}(\A))$ is an $e-$ultrafilter.\\
(d) Let $Z$ be a zero set that meets every
member of $E(E^{-}(\A))$ then, there exists a $\B\in [\A]$ such that
$Z\in \B$.
\end{lem}
{\flushleft\bf Proof.} (a) By Remark 2.2(v), pick $\mu\in\s$ such that $\bigcap_{A\in\A}\overline{\varepsilon(A)}=\{\mu\}$.
Now let $Z(f)\in\A$, then $\mu\in\overline{\varepsilon(Z(f))}$ and so there exists a net $\{\varepsilon(x_\alpha)\}\subseteq \varepsilon(A)$
such that $lim_\alpha \varepsilon(x_\alpha)=\mu$. Since
$$\mu(f)=lim_{\alpha}\varepsilon(x_\alpha)(f)=lim_\alpha f(x_\alpha)=0,$$
so $f\in Ker(\mu)$. It is obvious $Z(f)\subseteq \E(f)$ for each $\epsilon>0$ and so $\E(f)\in\A$ for each $\epsilon>0$.
This implies $f\in E^{-}(\A)$.

(b) By (a), there exists a $\mu\in\s$ such that $Ker(\mu)\subseteq E^{-}(\A)$. Since $Ker(\mu)$
is a  maximal ideal in $\LM(S)$ and also by Lemma 2.8(3), so \linebreak $Ker(\mu)=E^{-}(\A)$.

(c) Since $E^{-}(\A)$ is a maximal ideal so $E(E^{-}(\A))$ is an $e-$ultrafilter by Theorem 2.9.

(d) Let $Z$ be a zero set that meets every
member of $E(E^{-}(\A))$. Then $\{Z\}\cup E(E^{-}(\A))$ has
the finite intersection property. Hence there exists some $z-$ultrafilter $\B$ containing $\{Z\}\cup E(E^{-}(\A))$.
Since $E(E^{-}(\A))$ is an $e-$ultrafilter contained in $\B$, so by (b), $E^{-}(\B)$ is a maximal ideal and by Lemma 2.8(4),
$E^{-}(\A)\subseteq E^{-}(\B)$. Thus by Theorem 2.9, $E^{-}(\B)=E^{-}(\A)$ and so $E(E^{-}(\B))=E(E^{-}(\A))$. Therefore there exists a
$\B\in [\A]$ such that $Z\in\B$.$\hfill\square$

\vspace{.3cm}
\begin{rem}
Since $(\mathbb{R},+)$ is a locally compact topological group, by Theorem 5.7 of chapter 4  in \cite{Analyson},
$$\LM(\mathbb{R})=\{f\in \mathcal{CB}(\mathbb{R}):t\mapsto f\circ \lambda_t:\mathbb{R}\rightarrow \mathcal{CB}(\mathbb{R})\mbox{ is norm continuous.}\}.$$
 Let $\mathcal{C}_\circ(\mathbb{R})=\{f\in \mathcal{CB}(\mathbb{R}): lim_{x\rightarrow \pm\infty}f(x)=0\}$, then
 $\mathcal{C}_\circ(\mathbb{R})$ is an ideal of $\LM(\mathbb{R})$. Let $M$ be a maximal ideal in $\LM(\mathbb{R})$ which contains
 $\mathcal{C}_\circ(\mathbb{R})$. It is obvious that $f(x)=e^{-x^2}sin(x)$ and $g(x)=e^{-x^2}cos(\pi x)$ belong to $C_\circ(\mathbb{R})$. Then
 $Z(f)=\{k\pi:k\in \mathbb{Z}\}$, $Z(g)=\{\frac{2k+1}{2}:k\in\mathbb{Z}\}$, and $E(M)\cup\{Z(f)\}$
 and $E(M)\cup\{Z(g)\}$ have the finite intersection property. So there exist z-ultrafilters
 $\A$ and $\B$ such that $E(M)\cup\{Z(f)\}\subseteq \A$ and also \linebreak $E(M)\cup\{Z(g)\}\subseteq\B$.
 Since $E(M)$ is an $e-$ultrafilter so there exist at least two distinct $z-$ultrafilters
   containing $E^{-}(\A)$. Thus \\
   (i) If $\A$ is a $z-$ultrafilter, then it is not necessity the unique $z-$ultrafilter containing $E(E^{-}(\A))$.\\
   (ii) There exist a $z-$ultrafilter $\A$, and a zero-set $Z$ meets every member of $E(E^{-}(\A))$ such that
   $Z\notin\A$ .
\end{rem}

\section{\bf Space of e-ultrafilters}

In this section we define a topology on the set of all $e-$ultrafilters of a semitopological semigroup $S$,
and establish some of the properties of the resulting space. Also, we show that we can extend the operation
of the semitopological semigroup to the set of all $e-$ultrafilters.

\begin{defn}
Let $S$ be a Hausdorff semitopological semigroup.\\
$(a)$ The collection of all $e-$ultrafilters is denoted by $\mathcal{E}(S)$ i.e. $$\mathcal{E}(S)=\{p: p \mbox{ is an }e-\mbox{ultrafilter.}\}.$$
$(b)$  $A^{\dag}=\{p\in \mathcal{E}(S):A\in p\}$ for each $A\in Z(\LM(S))$ is defined.\\
$(c)$   $e(a)=\{\E (f):f(a)=0, \epsilon>0\}$ for each $a\in S$ is defined.\\
(d) It is said that $\A\subset Z(\LM(S))$ has the $e-$finite intersection property if and only if  $E(E^{-}(\A))$ has the finite
intersection property.
\end{defn}
Pick $\varepsilon(a)\in\s$ for some $a\in S$, then
\begin{eqnarray*}
Ker(\varepsilon(a))&=&\{f\in\LM(S):\varepsilon(a)(f)=0\}\\
&=&\{f\in\LM(S):f(a)=0\}
\end{eqnarray*}
is a maximal ideal and by Theorem 2.9,
$$ E^{-}(Ker(\varepsilon(a)))=\{\E(f):f(a)=0,\,\,\forall \epsilon>0\}=e(a)$$
 is an $e-$ultrafilter.
\begin{lem}
Let $A,B\in Z(\LM(S))$ and $f,g\in\LM(S)$. Then\\
$(1)$ ${(A\cap B)^\dag}={A}^\dag\cap{ B}^\dag$.\\
$(2)$ ${(A\cup B)}^\dag\supseteq{A}^\dag\cup{ B}^\dag$.\\
$(3)$ Pick $x\in S$ and $\epsilon>0$. Then
$\lambda_x^{-1}(\E(f))=\E(Lx f)$.\\
$(4)$  $E_{\epsilon\wedge\delta}(|f|\vee |g|)\subseteq\E(f)\cap
\Ed(g)$ and $\E(|f|\vee |g|)=\E(f)\cap \E(g)$, for each $\delta ,\epsilon>0$.
\end{lem}
{\flushleft\bf Proof.} The  proofs  are routine.$\hfill\square$
\vspace{.3cm}

 Since  ${(A\cap
B)^\dag}={A}^\dag\cap{ B}^\dag$, for each $A,B\in Z(\LM(S))$, so  the sets
${A}^\dag$ are closed under finite intersection.
Consequently, $\{{A}^\dag:A\in Z(\LM(S))\}$ forms a base for an open
topology on $\mathcal{E}(S)$.

\begin{thm}
(1) Pick $f\in \LM(S)$ and $\epsilon>0$, then $int_S(A)=e^{-1}(A^\dag)$, and so $e:S\rightarrow \mathcal{E}(S)$ is continuous.
\\
(2) Pick $p\in \mathcal{E}(S)$ and $A\in Z(\LM(S))$, then the following statements are equivalent:

(i) $p\in cl_{\mathcal{E}(S)}(e(A))$,

(ii) for each $B\in p$, $int_S(B)\cap A\neq\emptyset$,

(iii) for each $B\in p$, $B\cap A\neq\emptyset$,

(iv) there exists a $z-$ultrafilter $\A_p$ containing $p$ such that $A\in\A_p$.\\
(3) Pick $A,B\in Z(\LM(S))$ such that $p\in cl_{\mathcal{E}(S)}(e(A))\cap cl_{\mathcal{E}(S)}(e(B))$ and
$p\cup\{A,B\}$ has the finite intersection property, then $p\in
cl_{\mathcal{E}(S)}(e(A\cap B))$.
\\
(4)  $\{cl_{\mathcal{E}(S)}(e(A)): A\in Z(\LM(S))\}$ is a base for closed subsets of $\mathcal{E}(S)$.
\\
(5) $\mathcal{E}(S)$ is a compact Hausdorff space.
\\
(6) $e(S)$ is a dense subset
of $\mathcal{E}(S)$.
\end{thm}
{\flushleft\bf Proof.} (1) Let $p\in A^\dag$, so there exist $f\in E^{-}(p)$ and $\epsilon>0$ such that
 $E_\epsilon(f)=A$ and $E_\delta(f)\in p$ for each $\delta>0$. Pick $x_\circ\in int_S(A)$, then $|f(x_\circ)|<\epsilon$ or $|f(x_\circ)|=\epsilon$.

  If $\delta=|f(x_\circ)|<\epsilon$, then $E_{\epsilon-\delta}(|f|\vee \delta-\delta)=E_\epsilon(f)$,
 $x_\circ\in E_{\epsilon-\delta}(|f|\vee \delta-\delta)$ and $x_\circ\in E_{\eta}(|f|\vee \delta-\delta)$
 for each $\eta>0$. Thus,
 $$e(x_\circ)\in E_{\epsilon-\delta}(|f|\vee \delta-\delta)^\dag=E_\epsilon(f)^\dag=A^\dag.$$

 If $|f(x_\circ)|=\epsilon$ so there exists a neighborhood $U$ such that
 $x_\circ\in U\subseteq A$. Since $\LM(S)$ and $C(S^{\LM})$ are isometrically isomorphism,
 pick $g\in \LM(S)$ such that $g(U)=\{0\}$, $g(A^c)=\{\|f\|\}$ and $g(S)\subseteq [0,\|f\|]$. Define $h=|f|\wedge g$, then
 $\E(h)=\E(f)=A$ and $|h(x_\circ)|=0<\epsilon$. It is obvious that $E_\delta(f)\subseteq E_\delta(h)$ for each $0<\delta<\epsilon$ and
 $\E(f)\subseteq E_\delta(h)$ for each $\epsilon<\delta$. Therefore $E_\delta(h)\in p$ for each $\delta>0$ and $|h(x_\circ)|=0<\epsilon$.
 So by previous case, $e(x_\circ)\in E_\epsilon(h)^\dag=A^\dag$.
 Thus $x_\circ\in e^{-1}(A^\dag)$ and so $int_S(A)\subseteq e^{-1}(A^\dag)$.

 Now pick
$e(x)\in A^\dag$, so there exist $\epsilon>0$ and
$f\in\LM(S)$ such that $\E (f)=A$ and so, $E_{\delta}(f)\in e(x)$
for any $\delta>0$. Therefore, $f(x)=0$ and $x\in \E(f)$ for each $\epsilon>0$.
 Thus, $e^{-1}(A^\dag)=int_S(A)$.

(2) $(i)$ and $(ii)$ are equivalent. Since $p\in
cl_{\mathcal{E}(S)}(e(A))$ if and only if  $B^\dag\cap
e(A)\neq\emptyset$ for any $B\in p$, if and only if
$e^{-1}(B^\dag\cap e(A))\neq\emptyset$ for any $B\in p$, if and only if
$$int_S(B)\cap A=e^{-1}(B^\dag)\cap e^{-1}(e(A))\neq\emptyset$$
for any $B\in p$, by item (1).\\
It is obvious that $(iii)$ and $(iv)$ are equivalent and $(ii)$ implies $(iii)$.\\
$(iii)$ implies $(ii)$. Let for some $B\in p$, $B\cap A\neq\emptyset$ and $int_S(B)\cap A=\emptyset$.
Since $B\in p$ so there exist $f\in\LM(S)$ and $\epsilon>0$ such that $B=\E(f)$, $E_\delta(f)\in p$ for each $\delta>0$ and
 $$E_{\frac{\epsilon}{2}}(f)\cap A\subseteq int_S(B)\cap A=\emptyset.$$
This is a contradiction.\\

(3) Let $p\cup\{A,B\}$ has the finite intersection property, so $p\cup\{A\cap B\}$ has the finite intersection
property. Let $\A_p$ be a $z-$ultrafilter containing $p\cup\{A\cap B\}$ and hence by item (2), implies that $p\in
cl_{\mathcal{E}(S)}(e(A\cap B))$.

(4) It is suffices show that $\{(cl_{\mathcal{E}(S)}(e(A)))^c: A\in Z(\LM(S))\}$
is a base for open subsets of $\mathcal{E}(S)$. Let $U$ be an open subset
containing $p\in \mathcal{E}(S)$. Since $\{{A}^\dag:A\in Z(\LM(S))\}$ forms a base for an open
topology on $\mathcal{E}(S)$, so there exist $f\in \LM(S)$ and $\epsilon>0$ such that
$p\in\E(f)^\dag\subseteq U$ and $E_\delta(f)\in p$ for each $\delta>0$.
Now pick $0<\gamma<\mbox{min}\{\frac{\epsilon}{2}, \|f\|\}$, and define $g(x)=\|f\|-|f(x)|$. Then $g\in \LM(S)$ and
$(E_{\|f\|-\gamma}(g))^c\subseteq E_\gamma(f)$, so
$$(cl_{\mathcal{E}(S)}(E_{\|f\|-\gamma}(g)))^c\subseteq cl_{\mathcal{E}(S)}((E_{\|f\|-\gamma}(g))^c)\subseteq cl_{\mathcal{E}(S)}(E_\gamma(f)).$$
 Hence, there exists $\delta>0$ such that
$(E_{\|f\|-\gamma}(g)\cap E_{\gamma}(f))\bigcap E_\delta(f)=\emptyset$, and
$$E_{\|f\|-\gamma}(g)\bigcap E_\delta(f)=\emptyset.$$ This implies
$p\notin cl_{\mathcal{E}(S)}E_{\|f\|-\gamma}(g)$ and
so $$p\in (cl_{\mathcal{E}(S)}(E_{\|f\|-\gamma}(g))^c\subseteq \E(f)^\dag.$$
This show that $\{(cl_{\mathcal{E}(S)}(e(A)))^c: A\in Z(\LM(S))\}$
is a base for open subsets of $\mathcal{E}(S)$.

 (5) Suppose that $p$ and $q$ are distinct elements of $\mathcal{E}(S)$, then $E^{-}(p)$ and $E^{-}(q)$
 are maximal ideals, by Theorem 2.9. Pick $f\in E^{-}(p) \backslash E^{-}(q)$. So by Theorem 2.10, there exist $\epsilon>0$ and
 $A\in q=E(E^{-}(q))$  that $\E(f)\cap A=\emptyset$. Since $A\in q=E(E^{-}(q))$, pick
 $\delta>0$ and $g\in E^{-}(q)$ such that $A=E_{\delta}(g)$ and for all $\gamma>0$, $E_\gamma(g)\in q$. Then $\E(f)\cap E_\delta(g)=\emptyset$.
 Now let $B=E_\epsilon(f)$, then $A\in p$, $B\in q$ and $A\cap B=\emptyset$. Thus $A^\dag\cap
B^\dag=\emptyset$, $p\in{A}^\dag$ and $q\in{B}^\dag$, and so $\mathcal{E}(S)$ is Hausdorff.

Define $\eta: p\mapsto E(E^{-}(p)):\mathcal{Z}(S)\rightarrow \mathcal{E}(S)$. By Lemma 2.11,
if $p\in\mathcal{Z}(S)$ then $E(E^{-}(p))\in\mathcal{E}(S)$ so $\eta$ is well defined. Now let $p$
be an $e-$ultrafilter, so there exists a $z-$ultrafilter $\A$ containing $p$. By Lemma 2.11, \mbox{$p=E(E^{-}(\A))$}. This implies
$\eta$ is onto. For each $A\in
Z(Lmc(S))$, we have
\begin{eqnarray*}
\eta^{-1}(cl_{\mathcal{E}(S)}(e(A)))&=&\{p\in\mathcal{Z}(S):\eta(p)\in cl_{\mathcal{E}(S)}(e(A))\}\\
\mbox{By Theorem 3.3(2) }&=&\{\, p\in \mathcal{Z}(S):\forall B\in\eta(p),\,\,\,\,B\cap A\neq\emptyset\}\\
\mbox{By Theorem 3.3(2) }&=&\{\, p\in \mathcal{Z}(S):\eta(p)\cup\{A\}\subseteq p\}\\
 &=&\{\, p\in \mathcal{Z}(S):A\in p\}\\
 &=&\widehat{A}.
\end{eqnarray*}
Since $\{cl_{\mathcal{E}(S)}(e(A)): A\in Z(\LM(S))\}$ is a base for
closed subsets of $\mathcal{E}(S)$, so $\eta$ is continuous.
 Since $\mathcal{Z}(S)$ is compact, so $\mathcal{E}(S)$ is also  compact.

(6) By $(4)$, $e$ is continuous. Also,
\begin{eqnarray*}
\overline{e(S)}&=&\{p\in \mathcal{E}(S):\forall \,\, B\in p,\,\,\,
B^\dag\cap e(S)\neq\emptyset\}\\
&=&\{p\in e(S):\forall \,\, B\in p,\,\,\,
B\cap S\neq\emptyset\}\\
&=& \mathcal{E}(S).
\end{eqnarray*}
$ \hfill\square $
\begin{defn}
Let $\A$ be an $e-$filter. Then $\widehat{\A}=\{p\in
\mathcal{E}(S):\A\subseteq p\}$.
\end{defn}

\begin{thm}
$(a)$ If $\A$ is an $e-$filter, then $\widehat{\A}$ is a closed
subset of $\mathcal{E}(S)$.\\
$(b)$ Let $\A$ be an $e-$filter and $A\in Z(\LM(S))$. Then,
$A\in \A$ if and only if $\widehat{\A}\subseteq A^\dag$.\\
 $(c)$
Suppose that $A\subseteq \mathcal{E}(S)$ and $\A=E(E^{-}(\cap A))$, then $\A$ is an $e-$filter
and  $\widehat{\A}=cl_{\mathcal{E}(S)}A$.

\end{thm}
{\flushleft\bf Proof.} (a) Pick $p\in cl_{\mathcal{E}(S)}\widehat{\A}$, so $A^\dag\cap\widehat{\A}\neq\emptyset$, for each $A\in p$. Hence,
$\A\cup\{A\}$ has the e-finite intersection property for each $A\in p$. This implies that $\A\cup p\subseteq p$ and so $p\in\widehat{\A}$.

(b) It is routine to verify the assertion.

(c) By assumption, $\A$ is an $e-$filter (by Lemma 2.8). Further, for each $p\in A$, $\A\subseteq p$ implies that
$A\subseteq \widehat{\A}$, thus by $(a)$,
$cl_{\mathcal{E}(S)}A\subseteq \widehat{\A}$.

To see that $\widehat{\A}\subseteq cl_{\mathcal{E}(S)}A$, let $p\notin
cl_{\mathcal{E}(S)}A$. Then, there exist $B\in p$ and $C\in
Z(Lmc(S))$ such that $cl_{\mathcal{E}(S)}A\subseteq C^\dag$ and
$B^\dag\cap C^\dag=\emptyset$. Hence, $\widehat{\A}\subseteq
C^\dag$ and this implies $p\notin\widehat{\A}$.$\hfill\square$

\begin{defn}
Suppose that $p,q\in \mathcal{E}(S)$ and $A\in Z(\LM(S))$. Then, $A\in p+q$ if
there exist $\epsilon>0$ and $f\in\LM(S)$ such that
$A=\E(f)$ and $E_{\delta}(q,f)=\{x\in
S:\lambda^{-1}_{x}(E_{\delta}(f))\in q\}\in p$ for each $\delta>0$.
\end{defn}
\begin{thm}
Let $p,q\in \mathcal{E}(S)$, then $p+q$ is an $e-$ultrafilter.
\end{thm}
 {\flushleft\bf Proof.} It is obvious that $\emptyset\notin p+q$ and $S\in p+q$.
 Let $A\in p+q$, then there exist $\epsilon>0$ and $f\in\LM(S)$
such that $A=\E(f)$ and for each $\delta>0$, $E_{\delta}(q,f)=\{x\in
S:\lambda^{-1}_{x}(E_{\delta}(f))\in q\}\in p$. Thus, $\Ed(f)\in p+q$ for each $\delta>0$.
 Let $A,B\in p+q$, therefore, there exist $\delta,\epsilon>0$ and
 $f,g\in\LM(S)$ such that $A=\E(f)$ and $B=\Ed(g)$. So
 \begin{eqnarray*}
A\cap B&=&\E(f)\cap \Ed(g)\\
&\supseteq& E_{\epsilon\wedge\delta}(f)\cap E_{\epsilon\wedge\delta}(g)\\
&=&E_{\epsilon\wedge\delta}(|f|\vee |g|),
\end{eqnarray*}
and
\begin{eqnarray*}
E_\gamma (q,|f|\vee |g|)&=&\{x\in S:\lambda_{x}^{-1}(E_\gamma
(|f|\vee |g|))\in q\}\\
&=&\{x\in S:E_\gamma (|L_x f|\vee |L_x g|)\in q\}\\
&=&\{x\in S:E_\gamma (L_x f)\cap E_\gamma (L_x g)\in q\}\\
&=&E_\gamma (q, f)\cap E_\gamma (q, g).
\end{eqnarray*}
Since $E_\gamma (q, f), E_\gamma (q, g)\in p$ so $E_\gamma (q,
|f|\vee |g|)= E_\gamma(q, f)\cap E_\gamma (q, g)\in p$.
 Thus, $E_{\delta\wedge\epsilon} (|f|\vee |g|)\in p+q$ and so $A\cap B\in p+q$.

Now pick $A\in p+q$ and $B\in Z(\LM(S))$ such that $A\subseteq B$. So $A\in
p+q$ implies that there exist $\epsilon>0$ and $f\in \LM(S)$
such that $\E(f)=A$ and $\Ed(q,f)\in p$ for each $\delta>0$. For $B\in
Z(\LM(S))$ so, there exists a $g\in\LM(S)$ such that $Z(g)=B$. Now
define $u(x)=g(x)+\frac{\epsilon}{|f(x)|\vee\epsilon}$. Clearly, $h=\frac{u}{\|u\|}\in\LM(S)$,
$Z(g)=\E(fh)$ and $L_x f\in E^{-}(q)$ for each $x\in \Ed(q,f)$ and $\delta>0$.
 This implies $L_x fL_x h\in E^{-}(q)$ for
each $x\in \Ed (q,f)$, and so $E_\gamma
(L_xfL_xh)\in q$ for each $\gamma>0$. Thus, $E_{\delta}(q,f)\subseteq \Ed(q,fh)$ and
$\Ed(q,fh)\in p$ for each $\delta>0$, therefore, $Z(g)=\E(fh)\in
p+q$. So, $p+q$ is an $e-$filter.

Now, it is proved that $p+q$ is an $e-$ultrafilter. Let $E^{-}(p)=Ker(\mu)$
and $E^{-}(q)=Ker(\nu)$ for  $\mu,\nu\in \s$. It is claimed that
$E^{-}(p+q)=Ker(\mu\nu)$, thus $p+q$ is an $e-$ultrafilter.
Pick $f\in Ker(\mu\nu)$, so $T_\nu f\in Ker(\mu)$ and for each
$\epsilon>0$,
\begin{eqnarray*}
\E(T_\nu f)&=&\{x\in S:|T_\nu f(x)|\leq \epsilon\}\\
&=&\{x\in S:|\nu(L_x f)|\leq \epsilon\}\\
&=&\{x\in S:|\widehat{L_x f}(\nu)|\leq \epsilon\}\\
&\in& p.
\end{eqnarray*}
It is obvious that
 $\{t\in S:|\widehat{L_x f}(t)|\leq \epsilon\}=\{t\in S:|L_x f(t)|\leq \epsilon\}=\E(L_x
f).$
Pick $\epsilon>0$. For each $x\in E_{\frac{\epsilon}{2}}(T_\nu f)$,
 $E_{\frac{\epsilon}{2}}((|L_xf|\vee \frac{\epsilon}{2})-\frac{\epsilon}{2})\subseteq\E(L_xf),$
 and $E_{\frac{\epsilon}{2}}((|L_xf|\vee \frac{\epsilon}{2})-\frac{\epsilon}{2})\in E(Ker(\nu))=q$,  so
 $$\E(T_\nu f)\subseteq\{x\in S:\E(L_xf)\in q\}=\E(q,f).$$
Thus, $\E(f)\in p+q$ for each $\epsilon>0$, and this completes the proof.$\hfill\square$

\begin{thm}
$\mathcal{E}(S)$ and $\s$ are topologically isomorphic.
\end{thm}

 {\flushleft\bf Proof.} $M$ is a maximal ideal of $\LM(S)$ if and only if there
is a $\mu\in\s$ such that $Ker(\mu)=M$. Thus, $\gamma
:\mu\mapsto E(Ker(\mu)):\s\rightarrow \mathcal{E}(S)$ is well defined and surjective.
By Lemma 3.3(4), $\{cl_{\mathcal{E}(S)}(e(A)): A\in Z(\LM(S))\}$ is a base for
closed subsets of $\mathcal{E}(S)$, pick $A\in
Z(\LM(S))$ then
\begin{eqnarray*}
\gamma^{-1}(cl_{\mathcal{E}(S)}e(A))&=&\{\mu\in \s: E(Ker(\mu))\in cl_{\mathcal{E}(S)}e(A)\}\\
&=&\{\mu\in \s:\forall B\in E(Ker(\mu)),\,\,\,B^\dag\cap e(A)\neq\emptyset\}\\
&=&\{\mu\in \s:\forall f\in Ker(\mu ),\,\forall\,\delta >0,\,\,E_{\delta}(f)\cap A\neq\emptyset\} \\
&=&\{\mu\in \s:\forall f\in Ker(\mu ),\,\forall\delta >0,\,\exists x_\delta\in A\cap E_{\delta}(f)\}\\
&=& cl_{\s}(A).
\end{eqnarray*}
So $\gamma$ is continuous. Since, $\gamma :\s\rightarrow \mathcal{E}(S)$ is a surjective continuous function, and $\s$ a compact space therefore,
$\gamma$ is homeomorphism.
Now pick $\mu,\nu\in\s$, then
\begin{eqnarray*}
\gamma(\mu\nu)&=&E(Ker(\mu\nu))\\
\mbox{(see the proof of Theorem 3.7) }&=&E(Ker(\mu))+E(Ker(\nu))\\
&=&\gamma(\mu)+\gamma(\nu).
\end{eqnarray*}
Therefore, $\gamma$ is homomorphism and thus $\mathcal{E}(S)$ and $\s$ are topologically isomorphic.$\hfill\square$

\vspace{.2cm}
 By Theorem 3.8, $\s$ could be described  as a space
of $e-$ultrafilters, i.e. $\s=\{E(Ker(\mu)) :\mu\in\s \}$.
\begin{lem}
Let $A\in Z(\LM(S))$ and $x\in S$. Then $A\in e(x)+p$ if and only
if $\lambda_x^{-1}(A)\in p$.
\end{lem}
{\flushleft\bf Proof.} Pick $A\in e(x)+q$, so there exist
$\epsilon>0$ and $f\in\LM(S)$ such that $A=\E(f)$ and
$E_{\delta}(q,f)=\{t\in S:\lambda^{-1}_{t}(E_{\delta}(f))\in
q\}\in e(x)$ for each $\delta>0$ and $\lambda_x^{-1}(\Ed(f))\in q$
for each $\delta>0$. This implies $\lambda_x^{-1}(A)\in p$.\\
Conversely, let $\lambda_x^{-1}(A)\in p$, so there exist $\epsilon>0$ and $f\in\LM(S)$
such that $A=\E(f)$ and $\lambda_x^{-1}(A)\in p$. Thus $E_\delta(L_xf)=\lambda_x^{-1}(E_\delta(f))\in p$ for each $\delta>0$, and
$L_xf\in E^{-}(p)=Ker(\mu)$ for some $\mu\in S^{\LM}$. Clearly, $\mu(L_xf)=0$ and so $\varepsilon(x)\mu(f)=0$. This implies
$A\in E(Ker(\varepsilon(x)\mu))=e(x)+p$.$\hfill\square$
\begin{defn}
Let $\A$ and $\B$ be $e-$filters, and pick $A\in Z(\LM(S))$. Then
$A\in \A+\B$ if there exist $\epsilon>0$ and
$f\in\LM(S)$ such that $\E(f)=A$ and $\Ed(\B,f)=\{x\in
S:\lambda_x^{-1}(\Ed(f)\in \B)\}\in \A$ for each $\delta>0$.
\end{defn}
\begin{lem}
Let $\A$ and $\B$ are $e-$filters. Then $\A+\B$ is an $e-$filter.
\end{lem}
{\flushleft\bf Proof.} See Theorem 3.7.$\hfill\square$

\section{\bf Applications}

In this section, as an application, we consider to the semigroup \linebreak $S^*=\s  \setminus S$ and obtain some conditions characterizing
when $S^*$ is a left ideal of $\s$.

\begin{thm}
Pick  $p,q\in \mathcal{E}(S)$ and let $f\in \LM(S)$. Then $\E(f)\in p+q$ for
each $\epsilon>0$ if and only if for each $\epsilon>0$ there
exist $B_\epsilon\in p$ and an indexed family
$<C_{\epsilon,s}>_{s\in B_\epsilon}$ in $q$ such that $\bigcup
sC_{\epsilon,s}\subseteq \E(f)$.
\end{thm}
{\flushleft\bf Proof.} Let $\E(f)\in p+q$ for each $\epsilon>0$. Pick $\epsilon>0$, $x\in B_\epsilon=\E(q,f)$ and let
$C_{\epsilon, x}=\E(L_xf)=\lambda_{x}^{-1}(\E(f))$. For each $x\in B_\epsilon$, $C_{\epsilon, x}\in q$ and so
$\bigcup_{x\in B_\epsilon}xC_{\epsilon, x}\subseteq \E(f)$.

Conversely, by hypothesis for each $\epsilon>0$, there exist
$B_\epsilon\in p$ and an indexed family
$<C_{\epsilon,s}>_{s\in B_\epsilon}$ in $q$ such that
$\bigcup_{s\in B_\epsilon}sC_{\epsilon, s}\subseteq \E(f)$.
Then for each $s\in B_\epsilon$, $C_{\epsilon, s}\subseteq
\lambda_s^{-1}(\E(f))=\E(L_sf)$ and so $\E(L_sf)\in q$, for each
$s\in B_\epsilon$. Thus $B_\epsilon\subseteq\{t\in
S:\E(L_tf)\in q\}=\E(q,f)\in p$, and $\E(f)\in p+q$ for
each $\epsilon>0$.$\hfill\square$

\begin{thm}
Let $\A\subseteq Z(\LM(S))$ has the $e-$finite intersection
property. If for each $A\in E(E^{-}(\A))$ and $x\in A$, there exists
$B\in E(E^{-}(\A))$ such that $xB\subseteq A$, then
$\bigcap_{A\in E(E^{-}(\A))}\overline{\varepsilon(A)}$ is a subsemigroup of $\s$.
\end{thm}
{\flushleft\bf Proof.} Let $T=\bigcap_{A\in E(E^{-}(\A))}\overline{\varepsilon(A)}$. Since $E(E^{-}(\A))$
has the $e-$finite \linebreak intersection property, so $T\neq\emptyset$. Pick $p,q\in T$ and let
$A\in E(E^{-}(\A))$. Given $x\in A$, there is some
$B\in E(E^{-}(\A))$ such that $xB\subseteq A$. Therefore, there exist $f,g\in
\LM(S)$ such that $B=\Ed(g)$, $A=\E(f)$ and $E_\gamma(g),
E_\gamma(f)\in p\cap q$ for each $\gamma>0$, so $x\Ed(g)\subseteq
\E(f)$ and $\Ed(g)\subseteq \lambda_x^{-1}(\E(f))=\E(L_xf)$.
Since $B\in p\cap q$ thus $A\subseteq \{t\in S:\E(L_tf)\in q\}=\E(q,f)$, and $A=\E(f)\in p+q$.$\hfill\square$
\vspace{.3cm}

For semitopological semigroup $S$, let $S^*=\s \setminus \varepsilon(S)$.
\begin{defn}
$(a)$ $A\subseteq S$ is an unbounded set if $\overline{\varepsilon(A)}\cap
S^*\neq\emptyset$.\\
$(b)$ A sequence $\{x_n\}$ is an unbounded sequence if
$\overline{\varepsilon(\{x_n:n\in \mathbb{N}\})}\cap S^{*}\neq\emptyset$.
\end{defn}
\begin{lem}
Let $\{x_n\}$ and $\{y_n\}$ be unbounded sequences in $S$. Let
$p,q\in S^*$, $q\in \overline{\varepsilon(\{x_n:n\in \mathbb{N}\})}$ and $p\in
\overline{\varepsilon(\{y_n:n\in \mathbb{N}\})}$, then
$$p+q\in
\overline{\varepsilon(\{y_kx_n:k<n,\,\,k,n\in \mathbb{N}\})}.$$
\end{lem}
{\flushleft\bf Proof.} It is obvious that for each $A\in q$, $\varepsilon(\{x_n:
n\in \mathbb{N}\}) \cap {A}^\dag\neq \emptyset$ and for each $B\in
p$, $\varepsilon(\{y_n:n\in \mathbb{N}\})\cap {B}^\dag\neq\emptyset$. Now let
$C\in p+q$, then there exist $\epsilon>0$ and $f\in\LM(S)$ such
that $C=\E(f)$ and for each $\delta>0$, $\Ed(q,f)\in p$. Pick $\delta>0$ and let
$x\in \Ed(q,f)$, then
$$\varepsilon(\lambda_x^{-1}(\Ed(f))\cap \{x_n:n\in \mathbb{N}\})$$
and
$$\varepsilon(\Ed(q,f)\cap\{y_n:n\in \mathbb{N}\})$$
 are unbounded, by Lemma 3.3(4). Hence for each
$y_k\in\Ed(q,f)$, $$\varepsilon(\lambda_{y_k}^{-1}(\Ed(f))\cap \{x_n:n\in
\mathbb{N}\})$$
and so
$$\varepsilon(\{y_kx_n:k,n\in \mathbb{N}, \,\,
k<n\}\cap \Ed(f))$$
 are unbounded, by Lemma 3.3(4). This implies $\varepsilon(\{y_kx_n:k,n\in
\mathbb{N}\})\cap C^\dag\neq\emptyset$ and $p+q\in
\overline{\varepsilon(\{y_kx_n:k<n,\,\,k,n\in
\mathbb{N}\}})$.$\hfill\square$

\begin{thm}
Suppose that $S$ be a commutative semigroup, then $\s$ is not commutative if and
only if there exist unbounded sequences $\{x_n\}$ and $\{y_n\}$ such that
$$\overline{\varepsilon(\{x_ky_n:k<n,\,\,k,n\in
\mathbb{N}\})}\cap\overline{\varepsilon(\{y_kx_n:k<n,\,\,k,n\in
\mathbb{N}\})}=\emptyset.$$
\end{thm}
{\flushleft\bf Proof.} Necessity. Pick $p$ and $q$ in $S^*$ such that
$p+q\neq q+p$. Then, there exist $A\in p+q$ and $B\in q+p$ such that
$\overline{\varepsilon(A)}\cap \overline{\varepsilon(B)}=\emptyset$. So, there exist $\gamma, \epsilon>0$ and
$f,g\in \LM(S)$ such that $\E(f)=A$ and $E_\gamma(g)=B$. Pick
$0<\delta<\epsilon\wedge\gamma$, let $A_1=\Ed(q,f)$ and
$B_1=\Ed(p,g)$. Then, $A_1\in p$ and $B_1\in q$. Choose $x_1\in
A_1$ and $y_1\in B_1$. Inductively given $x_1, x_2,...,x_n$ and
$y_1, y_2,...,y_n$, choose $x_{n+1}$ and $y_{n+1}$ such that
$$\varepsilon(x_{n+1})\in
\varepsilon({A_1}^\dag\cap(\bigcap_{k=1}^{n}\lambda_{y_k}^{-1}(\Ed(g))))$$
 and
$$\varepsilon(y_{n+1})\in
\varepsilon({B_1}^\dag\cap(\bigcap_{k=1}^{n}\lambda_{y_k}^{-1}(\Ed(f)))).$$
 Then, $\{x_n\}$ and $\{y_n\}$ are unbounded sequences,
$$\varepsilon(\{y_kx_n:k,n\in \mathbb{N}, \,\, k<n\})\subseteq \varepsilon(A)$$
 and
$$\varepsilon(\{x_ky_n:k,n\in \mathbb{N}, \,\, k<n\})\subseteq \varepsilon(B).$$

Sufficiency. Now let there exist unbounded sequences $\{x_n\}$ and $\{y_n\}$
such that
$$\overline{\varepsilon(\{x_ky_n:k<n,\,\,k,n\in
\mathbb{N}\})}\cap\overline{\varepsilon(\{y_kx_n:k<n,\,\,k,n\in
\mathbb{N}\})}=\emptyset.$$
Pick $p\in\overline{\varepsilon(\{x_n:n\in
\mathbb{N}\})}\cap S^*$ and $q\in \overline{\varepsilon(\{y_n:n\in \mathbb{N}\})}\cap S^*$. Then
by Lemma 4.4,
$$q+p\in\overline{\varepsilon(\{y_kx_n:k<n,\,\,k,n\in
\mathbb{N}\})}$$
 and
 $$p+q\in\overline{\varepsilon(\{x_ky_n:k<n,\,\,k,n\in
\mathbb{N}\})}.$$
$\hfill\square$

\begin{defn}
A semitopological semigroup $S$ is a topologically weak left
cancellative if for all $u\in S$ there exists a
compact zero set $A$ such that $\varepsilon(u)\in A^\dag$ and $\lambda_v^{-1}(A)$ be a
compact set for each $v\in S$.
\end{defn}
\begin{thm}
(a) Let $S$ be a locally compact non-compact Hausdorff  semitopological semigroup and $S^*$ be a closed left ideal of $\s$, then
 $S$ is a topologically weak left cancellative.
\\
(b) Let $S$ be a locally compact non-compact Hausdorff  semitopological \linebreak semigroup and $S$ be a
topologically weak left cancellative,  then $S^*$ is a left ideal of $\s$.
\\
(c) Let $S$ be a locally compact non-compact Hausdorff semitopological \linebreak semigroup and $S^*$ be a closed subset of $\s$.
Then $S^*$ is a left ideal of $\s$ if and only if $S$ is a
topologically weak left cancellative.

\end{thm}
{\flushleft\bf Proof.} (a) Pick $x,y\in S$ such that for each compact zero set $A\in
Z(\LM(S))$, $\varepsilon(x)\in A^\dag$ and $B_A=\lambda_y^{-1}(A)$ is non-compact. Pick
$p_A\in S^*\cap \overline{\varepsilon(B_A)}$ so $\varepsilon(y)+p_A\in \varepsilon(A)$. Now let
 $$\mathcal{U}=\{A\in
Z(\LM(S)):\varepsilon(x)\in {A}^\dag\mbox{ and }A\mbox{ is compact}\},$$
 then
$\{p_A\}_{A\in \mathcal{U}}$ is a net, $\varepsilon(y)+p_A\rightarrow \varepsilon(x)$, and $\varepsilon(x)\in\overline{ S^*}=S^*$. So this is a contradiction.

(b) Since $S$ is non-compact so $S^*\neq\emptyset$. Pick $p\in S^*$, $q\in \s$ and let $q+p=\varepsilon(x)\in \varepsilon(S)$.
Let $A\in Z(\LM(S))$ be a compact set and $\varepsilon(x)\in
{A}^\dag$. Then $A\in q+p$ and there exist $f\in
\LM(S)$ and $\epsilon>0$ such that $\E(f)=A$ and
$\Ed(p,f)\in q$ for each $\delta>0$. Now pick $y\in \E(p,f)$ then $\lambda_y^{-1}(A)\in p$, so
$\lambda_y^{-1}(A)$ is not compact and this is a contradiction.\\
(c) It can easily be verified.$\hfill\square$

\begin{cor}
Let $G$ be a locally compact non compact Hausdorff topological group. Then $G^*$ is a left ideal of $G^\mathcal{LUC}$ if and only if $G$ is a
topologically weak left cancellative.
\end{cor}
{\flushleft\bf Proof.} Let $G$ be a locally compact non compact Hausdorff topological group, so $\varepsilon(G)$
is an open subset of $G^\mathcal{LUC}$, and $G^*$ is closed. Now by Theorem 4.7 proof is completed.$\hfill\square$

\begin{thm}
Let $S$ be a locally compact semitopological semigroup. The
following statements are equivalent.
\\
(a) $S^*$ is right ideal of $\s$.
\\
(b) Given any zero compact subset $A$ of $S$, any sequence $\{z_n\}$
in $S$, and any unbounded sequence $\{x_n\}$ in $S$, there exists a
$n<m$ in $\mathbb{N}$ such that $x_n\cdot z_m\notin A$.
\end{thm}
{\flushleft\bf Proof.} (a) implies (b). Suppose that $\{x_n\cdot z_m:n,m\in
\mathbb{N}\mbox{ and }n<m\}\subseteq A$. Pick $p\in\overline{\varepsilon(\{z_m:m\in\mathbb{N}\})}$ and
$q\in S^*\cap\overline{\varepsilon(\{x_n:n\in\mathbb{N}\})}$, which one can do,
since $\{x_n:n\in\mathbb{N}\}$ is unbounded. Thus
$q+p\in \overline{\varepsilon(A)}=\varepsilon(A)\subseteq \varepsilon(S)$, is a
contradiction.

(b) implies (a). Since $S^*\neq\emptyset$, pick $p\in\s$ and $q\in S^*$ such that $q+p=\varepsilon(a)\in \varepsilon(S)$, so there
exists a compact set $A\in Z(\LM(S))$ such that $\varepsilon(a)\in A^\dag$.
Hence there exist $\epsilon>0$ and $f\in \LM(S)$
such that $\E(f)=A$ and $\Ed (f)\in \varepsilon(a)$, for each
$\delta>0$. Then for each $1/n<\epsilon$,
$$E_{1/n}(p,f)=\{s\in S:\lambda_s^{-1}(E_{1/n}(f))\in p\}\in q,$$
choose an unbounded sequence $\{x_n\}$ such that $x_n\in
E_{1/n}(p,f)$. \linebreak Inductively choose a sequence $\{z_m\}$ in $S$ such
that for each $m\in \mathbb{N}$, \linebreak
$z_m\in\bigcap_{n=1}^m\lambda_{x_n}^{-1}(E_{1/n}(f))$ (which one can
do) since $\bigcap_{n=1}^m\lambda_{x_n}^{-1}(E_{1/n}(f))\in p$. Then
for each $n<m$ in $\mathbb{N}$, $x_n\cdot z_m\in E_{1/n}(f)\subseteq
\E(f)=A$, is a contradiction.

 $\hfill\square$

{\flushleft\bf Examples.}
 (a) Let $S$ be a discrete semigroup. If $S$ is either right or left cancellative then $S^*=\beta S \setminus S$ is a
subsemigroup of $\beta S$, (See Corollary 4.29 in \cite{hindbook}). This is not true for a left cancellative semitopological semigroup $S$.
Let $(S=(1,+\infty),+)$  with the natural topology.
Then $S^*$ is not subsemigroup. Pick $p,q\in cl_{\s}(1,2]$, thus there exist nets $\{x_\alpha\}$ and $\{y_\beta\}$ in $(1,2]$
such that $x_\alpha \rightarrow p$, $y_\beta \rightarrow q$
and $x_\alpha +y_\beta\in [2,4]$. Hence $p+q\in [2,4]$ and so $S^*$ is not subsemigroup.
Also, $S^*$ is not a left ideal and so $S$ is not topologically weak left cancellative.

(b) $(S=[1,+\infty),+)$ with the natural topology is a topologically weak left cancellative, thus $S^*$ is a left ideal of $\s$.\\
\\
\\

\bibliographystyle{alpha}

\end{document}